\documentclass[a4paper]{amsart}
\usepackage[british]{babel}
\usepackage{amsmath}
\usepackage{amsthm}
\usepackage{amscd}
\usepackage{amsfonts}
\usepackage{amssymb}
\usepackage{hyperref}
\usepackage{ifthen}

\newcommand{\Q}{\mathbb{Q}}
\newcommand{\Z}{\mathbb{Z}}
\newcommand{\X}{\mathcal{X}}
\newcommand{\Xs}{\mathcal{X}_s}
\newcommand{\G}{\mathbb{G}}
\newcommand{\ad}{\mathbb{A}}
\renewcommand{\O}{\mathcal{O}}
\renewcommand{\P}{\mathbb{P}}
\newcommand{\knr}{{k^\textrm{nr}}}
\newcommand{\F}{\mathbb{F}}
\DeclareMathOperator{\inv}{inv}
\DeclareMathOperator{\Br}{Br}
\DeclareMathOperator{\res}{res}
\DeclareMathOperator{\Pic}{Pic}
\DeclareMathOperator{\Div}{Div}
\DeclareMathOperator{\Spec}{Spec}
\DeclareMathOperator{\Gal}{Gal}
\DeclareMathOperator{\Frac}{Frac}

\DeclareMathOperator{\Cl}{Cl}
\DeclareMathOperator{\Hom}{Hom}
\DeclareMathOperator{\ch}{char}

\newcommand{\kb}{\bar{k}}
\newcommand{\Xb}{\bar{X}}
\newcommand{\Xbs}{\bar{\mathcal{X}}_s}

\newcommand{\Xnr}{X^\textrm{nr}}
\newcommand{\Cb}{\bar{C}}
\newcommand{\Cbs}{\bar{\mathcal{C}}_s}
\newcommand{\Cs}{\mathcal{C}_s}

\newcommand{\al}[1]{\mathcal{#1}}
\newcommand{\A}{\al{A}}
\newcommand{\Y}{\mathcal{Y}}
\newcommand{\Ys}{\mathcal{Y}_s}
\newcommand{\Ybs}{\bar{\mathcal{Y}}_s}
\newcommand{\m}{\mathfrak{m}}

\renewcommand{\H}{\mathrm{H}}
\newcommand{\places}[1]{\mathfrak{M}_{#1}}

\newtheorem{theorem}{Theorem}
\newtheorem{lemma}[theorem]{Lemma}
\newtheorem{proposition}[theorem]{Proposition}
\newtheorem{corollary}[theorem]{Corollary}
\numberwithin{theorem}{section}
\theoremstyle{remark}
\newtheorem*{remark}{Remark}

\title{Evaluating Azumaya algebras on cubic surfaces} 

\author{Martin Bright}
\address{Mathematics Institute \\ Zeeman Building \\ University of Warwick \\
Coventry CV4 7AL \\ UK}
\email{M.Bright@warwick.ac.uk}
\thanks{This research was funded by the Heilbronn Institute for 
Mathematical Research.}


\begin{document}

\begin{abstract}
Let $X$ be a cubic surface over a $p$-adic field $k$.  Given an
Azumaya algebra on $X$, we describe the local evaluation map $X(k) \to
\Q/\Z$ in two cases, showing a sharp dependence on the geometry of the
reduction of $X$.  When $X$ has good reduction, then the evaluation
map is constant.  When the reduction of $X$ is a cone over a smooth
cubic curve, then generically the evaluation map takes as many values
as possible.  We show that such a cubic surface defined over a number
field has no Brauer--Manin obstruction.  This extends results of
Colliot-Th\'el\`ene, Kanevsky and Sansuc.
\end{abstract}

\maketitle

\section{Introduction}

Let $L$ be a number field, and let $X \subset \P^3_L$ be a smooth
cubic surface defined over $L$.  It is known that $X$ does not have to
satisfy the Hasse principle: that is, it is possible for $X(L_v)$ to
be non-empty for each place $v$ of $L$, but for $X$ nonetheless to
have no $L$-rational points; the first example of such a surface was
given by Swinnerton-Dyer~\cite{SD:M-1962}.  This and other
counterexamples to the Hasse principle were shown by
Manin~\cite{Manin:CF} to be explained by what is now known as the
Brauer--Manin obstruction.  It has been conjectured by
Colliot-Th\'el\`ene that the Brauer--Manin obstruction is in fact the
only obstruction to the Hasse principle for cubic surfaces (and, more
generally, for geometrically rational varieties:
see~\cite[p.~319]{CT:AFST-1992}).

Let $\places{L}$ denote the set of places of $L$.  Let $\Br X$ denote
the Brauer group of $X$, and $X(\ad_L)$ the set of adelic points;
since $X$ is projective, $X(\ad_L)$ is the same as $\prod_{v \in
  \places{L}} X(L_v)$.  The Brauer--Manin obstruction is based on the
pairing
\begin{equation} \label{eq:brpairing}
X(\ad_L) \times \Br X \to \Q/\Z, \qquad ((x_v)_{v \in \places{L}},
\al{A}) \mapsto \sum_{v \in \places{L}} \inv_v \al{A}(x_v) \text{.}
\end{equation}
Here $\inv_v \colon \Br L_v \to \Q/\Z$ is the invariant map of local
class field theory: see~\cite[XIII, \S 3]{Serre:LF}.  To understand
the obstruction, then, it is desirable to have a good description of
the local evaluation map $X(L_v) \to \Q/\Z$, $x \mapsto \inv_v \A(x)$
given by a particular element $\A$ of the Brauer group of $X$ at some
given place $v$.  For a general variety, not much is known about this
map, except that it is constant for $v$ outside some finite set of
primes; the usual way to compute it is simply to list the points of
$X(L_v)$ to sufficient accuracy and evaluate the invariant at each
one.  However, for curves the picture is much clearer: if $C$ is a
curve over a $p$-adic field $k$, then the local pairing extends to a
pairing $\Pic C \times \Br C \to \Br k$, which
Lichtenbaum~\cite{Lichtenbaum:IM-1969} showed can be identified with
that arising from the Tate pairing, and is therefore non-degenerate.

In~\cite{CTKS}, Colliot-Th\'el\`ene, Kanevsky and Sansuc gave a thorough
description of both $\Br X$ and the local evaluation maps in the case
when $X$ is a \emph{diagonal} cubic surface.  In particular, they
showed that the evaluation map is constant for places $v$ where $X$
has good reduction or where $X$ is rational over $L_v$; and, in
contrast (Proposition~2), at finite places $v$ where the reduction of
$X$ at $v$ is a cone, the local evaluation map takes all possible
values on $X(L_v)$.  In this last case the Brauer--Manin obstruction
therefore vanishes.  The proof relies on explicitly determining an
Azumaya algebra which generates $\Br X / \Br L$, showing that it has
non-trivial restriction to a certain nonsingular cubic curve, and
applying Lichtenbaum's result.

The object of this article is to extend some of the results
of~\cite{CTKS} to more general cubic surfaces $X$, and to show how
these results can be deduced from the geometry of a model of $X$.  In
particular, it is still straightforward to see that the local
evaluation map is constant at places of good reduction: although this
is already known, it gives a useful illustration of our approach; we
prove this in Theorem~\ref{thm:good}.  At places where $X$ reduces to
a cone, the behaviour described in~\cite{CTKS} still happens as long
as the singularity of the model of $X$ is not too severe; this is
proved in Theorem~\ref{thm:cone}, which relies on a description of the
geometry of the model given in Proposition~\ref{prop:flexes}.

\subsection{Background}

We now review the definition of the Brauer--Manin obstruction, partly
in order to fix notation.  An excellent reference for this topic is
Skorobogatov's book~\cite{Skorobogatov:TRP}.

Define the Brauer group of a scheme $X$ to be the \'etale cohomology
group $\Br X = \H^2(X, \G_m)$; for a smooth variety $X$, this is equal
to the group of equivalence classes of Azumaya algebras on $X$.  If
$K$ is any field, then a $K$-point of $X$ corresponds to a morphism
$\Spec K \to X$ and so by functoriality gives a homomorphism $\Br X
\to \Br K$.  In this way we obtain an ``evaluation'' map $X(K) \times
\Br X \to \Br K$.  In particular, suppose that $X$ is a variety over a
number field $L$; then, for each place $v$ of $L$, there is a map
$X(L_v) \times \Br X \to \Br L_v$, which we may compose with the local
invariant map $\inv_v \colon \Br L_v \to \Q/\Z$.  Let $X(\ad_L)$
denote the set of adelic points of $X$, which is equal to the product
$\prod_v X(L_v)$ if $X$ is projective.  Adding together the local
evaluation maps gives the map~\eqref{eq:brpairing}.
Manin~\cite{Manin:GBG} observed that the $L$-rational points of $X$
must lie in the left kernel of this map, and that this explained many
known counterexamples to the Hasse principle -- that is, varieties $X$
with $X(L_v) \neq \emptyset$ for all $v$, yet $X(L)=\emptyset$.  This
obstruction is known as the Brauer--Manin obstruction to the Hasse
principle.

If $X$ is any smooth, proper, geometrically integral variety over
any field $K$, there is an exact sequence as follows, which arises from
the
Hochschild--Serre spectral sequence for $\G_{m,X}$:
\begin{equation} \label{eq:hs}
\Br K \to \ker(\Br X \to \Br \Xb) \xrightarrow{r} \H^1(K, \Pic \Xb) \to
\H^3(K,\bar{K}^\times) \text{.}
\end{equation}

Here $\bar{K}$ denotes a fixed algebraic closure of $K$, and $\Xb$ is
the base change of $X$ to $\bar{K}$.  When $K$ is a number field or a
local field, we have $\H^3(K, \bar{K}^\times)=0$ and the homomorphism
$r$ is surjective.  We write $\Br_1 X$ for $\ker( \Br X \to \Br \Xb)$;
the map $r \colon \Br_1 X / \Br K \to \H^1(K, \Pic \Xb)$ is an
isomorphism.  If $\Xb$ is a rational variety, then $\Br \Xb=0$
(see~\cite[Theorem~42.8]{Manin:CF}) and therefore $\Br_1 X = \Br X$.
In this case $\H^1(K, \Pic\Xb)$ is finite and contains all the
interesting information about the Brauer group of $X$.

The Hochschild--Serre spectral sequence is functorial.  If particular,
if $Y \to X$ is a morphism of varieties, then there are natural maps
$\Br X \to \Br Y$ and $\Pic \Xb \to \Pic \bar{Y}$, and the following
diagram commutes:
\[
\begin{CD}
\Br K @>>> \Br_1 X @>>> \H^1(K, \Pic \Xb)\\
@|         @VVV         @VVV \\
\Br K @>>> \Br_1 Y @>>> \H^1(K, \Pic \bar{Y})
\end{CD} \text{.}
\]

When $X$ is a smooth cubic surface, it is well known that there are
exactly 27 straight lines on $X$, and that their classes generate the
Picard group of $\Xb$, which is isomorphic to $\Z^7$.
Swinnerton-Dyer~\cite{SD:MPCPS-1993} has described the structure of
$\H^1(K, \Pic\Xb)$ for all possible Galois actions on the 27 lines,
thus giving the structure of $\Br X / \Br K$ in each case.  It should
be noted that turning an explicit element of $\H^1(K, \Pic\Xb)$ into
an explicit element of $\Br X$ is in general a highly non-trivial
procedure; fortunately we will not have to do it.

\subsection{Notation}

Throughout, $k$ will be a finite extension of $\Q_p$.  We denote by
$\O$ the ring of integers of $k$, which has maximal ideal $\m$
generated by a uniformising element $\pi$.  The residue field $\O/\m$
is denoted by $\F$.  Let $\kb$ be a fixed algebraic closure of $k$ and
$\knr$ the maximal unramified extension of $k$ in $\kb$.  If $X$ is a
variety over $k$, we denote by $\Xnr$ and $\Xb$ the base changes of
$X$ to $\knr$ and $\kb$ respectively.  If $Y$ is a variety over $\F$,
then $\bar{Y}$ denotes the base change of $Y$ to the algebraic closure
of $\F$.  If $\al{A}$ is an element of $\Br X$, we denote the
associated evaluation map also by $\al{A} \colon X(k) \to \Q/\Z$.

\subsection{Reduction of projective varieties}
\label{sec:reduction}

Given a closed subvariety $X$ of $\P^3_k$, the \emph{reduction} $\Xs$
of $X$ is defined as follows: let $\X$ be the Zariski closure of $X$ in
$\P^3_\O$; then $\Xs$ is the special fibre of $\X$.  We can
characterise $\X$ as the unique closed subvariety of $\P^3_\O$ which
has generic fibre $X$ and is flat over $\O$ (see~\cite[III,
  Proposition~9.8]{Hartshorne:AG}).  If $X$ is a hypersurface defined
by a homogeneous polynomial $f$, then $\X$ is found simply by
multiplying $f$ by an appropriate power of $\pi$ so that the
coefficients of $f$ lie in $\O$ but are not all in $\m$, and $\X_s$ by
reducing the resulting polynomial modulo $\m$.

if $X$ and $Y$ are two closed subvarieties of $\P^3_k$, denote their
closures in $\P^3_\O$ by $\X$ and $\mathcal{Y}$ respectively.  Then
the scheme-theoretic intersection $(\X \cap \mathcal{Y})$ is closed in
$\P^3_\O$ and contains $X \cap Y$, and so contains the closure of $X
\cap Y$.  Therefore $(\X_s \cap \mathcal{Y}_s)$ contains the reduction
of $X \cap Y$.  This inclusion can be strict when $\X \cap
\mathcal{Y}$ is not flat over $\O$: consider, for example, the case
when $X$ and $Y$ are two lines in $\P^2_k$ which have the same
reduction.

The process of taking reductions commutes with base change: let $\O
\to R$ be a local morphism of discrete valuation rings, and
$L=\Frac(R)$.  Then, since $R$ is flat over $\O$, $\X \times_\O R
\subset \P^3_R$ is flat over $R$ and has generic fibre $X \times_k L$,
so is equal to the closure of $X \times_k L$ in $\P^3_R$.  In
particular, this means that we can talk about the reduction of a
closed subvariety of $\P^3_{\kb}$ without worrying about which finite
extension of $k$ it is defined over.

Applying reduction to Weil divisors on $X$ gives a map of divisor
class groups $\Cl X \to \Cl \Xs$: see~\cite[p.~399]{Fulton:IT}.  The
reduction of a Cartier divisor is not necessarily again Cartier; but
this is true if $\X$ is smooth, and so in that case there is a
reduction map $\Pic X \to \Pic \Xs$.  This is the same as the
composition of natural maps on Picard groups $\Pic X \cong \Pic \X \to
\Pic Xs$, given by extending a line bundle on $X$ to $\X$ (possible
because the special fibre is a principal prime divisor) and then
restricting to $\Xs$.  Passing to the limit over finite extensions of
$k$, we obtain a Galois-equivariant reduction map $\Pic \Xb \to \Pic
\Xbs$.

If $X$ is a \emph{diagonal} cubic surface, and $\ch \F > 3$, then the
reduction $\Xs$ is either smooth; a cone over a smooth plane cubic
curve; a union of three planes; or a triple plane, depending on the
number of coefficients which lie in the maximal ideal $\m$.  In
general, the reduction of a cubic surface can have other types of
singularity: see~\cite{BW:JLMS-1979} for a comprehensive description
of singularities of cubic surfaces.

\section{Good reduction}

We first treat the places of good reduction -- that is, where the
reduction of the cubic surface $X$, in the sense of
Section~\ref{sec:reduction}, is again a smooth cubic surface. In this
case the situation is very simple.  It is well known that, given a
variety $X/k$ with a smooth, proper model $\X/R$, and an Azumaya
algebra $\al{A} \in \Br X$ also with ``good reduction'' in the sense
that $\al{A}$ extends to an element of $\Br \X$, we have $\al{A}(P)=0$
for all $P \in X(k)$: for the evaluation map at any point factors
through the trivial group $\Br R$. The following theorem shows that,
in our situation, the condition that $\al{A}$ have good reduction is
unnecessary.  Of course, we can then no longer expect the evaluation
map to be zero, for this is not true for a non-zero constant algebra;
but the map is constant.

I thank an anonymous referee for pointing out that this theorem could
be made significantly more general than its original form.

\begin{theorem} \label{thm:good}
Let $\X$ be a smooth, proper scheme over $R$; assume that the generic
fibre $X$ and the special fibre $\Xs$ are both geometrically integral.
Suppose that $\Pic \Xb$ is torsion-free, and that $\H^1(\Xbs,
\O_{\Xbs})=0$.  Then, for any $\al{A} \in \Br_1 X$, the associated
evaluation map $\al{A}\colon X(k) \to \Q/\Z$ is constant.
\end{theorem}

In particular, this applies whenever $\Xb$ and $\Xbs$ are rational
varieties, and so to our case of a cubic surface with good reduction.

\begin{remark}
The conclusion of this theorem can also be reached by other means in
some situations: for example, when $X$ is a rational surface, it
follows from the fact that the Chow group of 0-cycles on $X$ is
trivial; see~\cite{CT:IM-1983}.  Let $p = \ch \F$; for more general
$X$, the result can also be obtained from the purity theorem for the
Brauer group, proved by Gabber (see~\cite{Fujiwara:PC}), as long as
the order of $\al{A}$ in $\Br X$ is coprime to $p$: one deduces that
such $\al{A}$ can be extended to $\Br \X$ after adding a constant
algebra.  Our view is that Azumaya algebras split by an extension of
the base field are much simpler to understand than arbitrary Azumaya
algebras, and so deep results such as the purity theorem in all its
generality are not necessary in this case.  The benefit is that our
result also holds for the $p$-part of the Brauer group.
\end{remark}

Theorem~\ref{thm:good} follows from the following lemmas.

\begin{lemma} \label{lem:unram}
Under the conditions of Theorem~\ref{thm:good}, the Galois module
$\Pic \Xb$ is unramified, meaning that the inertia subgroup of
$\Gal(\kb/k)$ acts trivially on $\Pic \Xb$.
\end{lemma}
\begin{proof}
Since $\X$ is smooth, there is a Galois-equivariant reduction map
$\Pic \Xb \to \Pic \Xbs$; the fact that $\H^1(\Xbs, \O_{\Xbs})=0$
implies that this map is injective, by~\cite[Corollaire~3 to
  Th\'eor\`eme~7]{Grothendieck:GFGA}.  But the inertia subgroup of
$\Gal(\kb/k)$ acts trivially on $\Pic \Xbs$, and so acts trivially on
$\Pic \Xb$ as well.
\end{proof}

\begin{corollary}\label{cor:unram}
Under the conditions of Theorem~\ref{thm:good}, $\H^1(\knr, \Pic\Xb)=0$.
\end{corollary}
\begin{proof}
The action (which is that of the inertia group) is trivial and
$\Pic\Xb$ is torsion-free.
\end{proof}

Recall the map $r\colon \Br_1 X \to \H^1(k, \Pic\Xb)$ from
~\eqref{eq:hs}.  We say that a class $\alpha \in \H^1(k, \Pic\Xb)$ is
\emph{split} by a finite extension $k'/k$ if the image of $\alpha$
under the restriction map $\res: \H^1(k, \Pic\Xb) \to \H^1(k',
\Pic\Xb)$ is zero.

\begin{lemma}\label{lem:constant}
Let $\X$ be a smooth, proper scheme over $R$; assume that the generic
fibre $X$ and the special fibre $\Xs$ are both geometrically integral.
Let $\al{A} \in \Br_1 X$ be such that $r(\al{A}) \in \H^1(k, \Pic\Xb)$
is split by $\knr$.  Then the evaluation map $\al{A} \colon X(k) \to
\Q/\Z$ is constant.
\end{lemma}
\begin{proof}
We first show that $\al{A}$ itself is split by $\knr$.  Consider the
commutative diagram
\begin{equation}
\begin{CD}
\Br k @>>> \Br_1 X @>r>> \H^1(k, \Pic\Xb) \\
@VVV        @VVV         @VVV \\
\Br \knr @>>> \Br_1 \Xnr @>>> \H^1(\knr, \Pic\Xb)
\end{CD}
\end{equation}
where the vertical maps are all given by restriction, and the rows
come from~\eqref{eq:hs}.  Since $r(\al{A})$ restricts to $0$ in
$\H^1(\knr,\Pic\Xb)$ and $\Br \knr=0$, the algebra $\al{A}$ restricts to
$0$ in $\Br_1 \Xnr$ -- that is, $\al{A}$ is split by base extension to
$\knr$.

The remaining argument is a special case of Theorem~1
of~\cite{Bright:MPCPS-2007}, which we briefly summarise.  Let $K/k$ be
an unramified extension splitting $\al{A}$, and write $K(X)$ for the
function field of $X_K$.  Recall that, as described in~\cite[X, \S
  5]{Serre:LF}, to any class in $\Br k(X)$ split by the extension
$K(X)/k(X)$ we can associate a cocycle class in $\H^2(K/k,
K(X)^\times)$.  (Here we implicitly identify $\Gal(K(X)/k(X))$ with
$\Gal(K/k)$.)  In particular, this applies to $\al{A}$ (strictly, its
restriction to the generic point of $X$).  As $X$ is smooth, the fact
that $\al{A}$ is Azumaya implies moreover that we can represent
$\al{A}$ by a cocycle class
\[ 
\alpha \in \ker\left( \H^2(K/k, K(X)^\times) \to \H^2(K/k, \Div X_K) \right)
\text{.}
\]
Let $\X_K$ denote the base change of $\X$ to the ring of integers
$\O_K$ of $K$.  Since $\X_K$ is smooth over $\O_K$, every Weil divisor
on $\X_K$ is Cartier, and we have $\Div \X_K = \Div X_K \oplus \Z D$,
where $D$ is the prime divisor corresponding to the special fibre.
The image of $\alpha$ under the natural map $\H^2(K/k,
K(X)^\times) \to \H^2(K/k, \Div \X_K)$ therefore lies in
$\H^2(K/k, \Z D)$.  Since $\pi$ is also a uniformiser in $K$, and the
constant function $\pi$ cuts out the divisor $D$ on $\X_K$, there is a
commutative diagram
\[
\begin{CD}
K^\times        @>{v_K}>>  \Z \\
@VVV                      @VV{1 \mapsto D}V \\
K(X)^\times @>{\textrm{div}}>> \Div \X_K
\end{CD} \qquad \text{.}
\]
Taking cohomology gives
\[
\begin{CD}
\Br(K/k)                  @>>> \H^2(K/k, \Z) \\
@VVV                       @VV{1 \mapsto D}V \\
\H^2(K/k, K(X)^\times) @>>> \H^2(K/k, \Div \X_K)
\end{CD} \text{.}
\]
Now, since $K/k$ is unramified, the map $\Br(K/k) \to \H^2(K/k,\Z)$ is
an isomorphism.  So we may choose a constant algebra $\al{B} \in
\Br(K/k)$ such that $\al{A'} := \al{A} - \al{B}$ maps to $0$ in
$\H^2(K/k, \Div \X_K)$.

Now let $P$ be a point of $X(k)$; then, as $\X$ is proper, $P$ reduces
to a point $Q \in \X(\F)$.  Let $S = \Spec \O_{\X,Q}$, let $S_K$ be
its base change to $\O_K$, and write $A$ for the affine coordinate
ring of $S_K$; then the field of fractions of $A$ is $K(X)$.  Since
$\Pic S_K=0$, there is an exact sequence
\[
\H^2(K/k, A^\times) \to \H^2(K/k, K(X)^\times) \to \H^2(K/k, \Div S_K) \text{.}
\]
The divisor map $K(X)^\times \to \Div S_K$ is the composite of the
divisor map $K(X)^\times \to \Div \X_K$ with the natural map $\Div
\X_K \to \Div S_K$ given by the inclusion $S_K \subset \X_K$, and so
the image of $\al{A'}$ is zero in $\Div S_K$.  We conclude that
$\al{A}'$ may be represented by a cocycle in $\H^2(K/k, K(X)^\times)$
taking values in $A^\times$.  It follows that the evaluation
$\al{A}'(P)$ lies in the subgroup $\H^2(K/k, \O_K^\times) \subseteq
\Br(K/k)$, but since $K/k$ is unramified we have $\H^2(K/k,
\O_K^\times)=0$.  Therefore $\al{A}'(P)=0$, and so $\al{A}(P) =
\al{B}(P)$ which is the same for all $P \in X(k)$.
\end{proof}

\section{Reduction to a cone}

In this section we assume that the characteristic of $\F$ is not $3$.

Suppose that $\X = \{ F=0 \} \subset \P^3_\O$ is a cubic surface with
smooth generic fibre $X$, and special fibre a cone over a smooth cubic
curve.  By choosing coordinates such that the vertex of the cone lies
at $(0:0:0:1) \in \P^3_\F$, the equation of $\X$ becomes of the form
\begin{equation} \label{eq:cone}
F = f(X_0,X_1,X_2) + \pi^s g(X_0,X_1,X_2,X_3)
\end{equation}
with $f,g$ homogeneous polynomials of degree $3$ with coefficients in
$\O$, the reduction modulo $\m$ of $f$ defining a smooth cubic curve
over $\F$, and where $s>0$ is chosen maximally: that is, so that some
monomial in $g$ involving $X_3$ has its coefficient in $\O^\times$.
The results in this section will apply in the case that $g(0,0,0,1)$
is a unit in $\O$, or equivalently that the coefficient of $X_3^3$ in
$g$ is a unit; this condition should be understood as saying that the
singularity of $\X$ is not ``too bad''.  Another equivalent statement
is that the tangent cone to the total space $\X$ at the closed point
corresponding to $(0:0:0:1)$ is defined by $\pi^s$.  In particular,
this is true if $\X$ is diagonal, or if $\X$ is regular (in which case
$s=1$).

In this situation, if $s \ge 3$, then we can perform the change of
variables $X_i \mapsto \pi X_i \quad (i=0,1,2)$ and then remove the
resulting factor of $\pi^3$ to give a new model of $X$; this
corresponds to the geometric operation of blowing up $\X$ at the
closed point with ideal $(X_0,X_1,X_2,\pi)$ and then blowing down the
strict transform of the old special fibre, although we will not use
this description.  In this way we can always find a model of the
form~\eqref{eq:cone} with $s < 3$.  If we reach $s=0$, then $X$ has
good reduction:

\begin{lemma} \label{lem:good}
Suppose that $\X$ has an equation of the form~\eqref{eq:cone}, with
$g(0,0,0,1)$ a unit in $\O$ and $s \ge 1$, and suppose that $s$ is
divisible by 3.  Let $X$ denote the generic fibre of $\X$.  Then there
is a smooth model $\X' \subset \P^3_\O$ with generic fibre isomorphic
to $X$; in other words, $X$ has good reduction.
\end{lemma}
\begin{proof}
First suppose that $s=3$.  Performing the change of variables
described above leads to a model $\X'$ with equation
\[
\pi^{-3} F(\pi X_0, \pi X_1, \pi X_2, X_3) = 
f(X_0, X_1, X_2) + a X_3^3 + \pi h(X_0, X_1, X_2, X_3)
\]
where $a$ is the coefficient of $X^3$ in $g$, and $h$ is a new cubic
form with coefficients in $\O$.  The special fibre of $\X'$ is a
triple cover of $\P^2_\F$ branched over the smooth cubic curve $\{ f=0
\}$, and is easily checked to be smooth.

If $s>3$, then repeating this process several times gives the same
result.
\end{proof}

It follows from Proposition~\ref{prop:flexes} below that, inversely, if
there is a model with $s$ not divisible by $3$, then $X$ never has
good reduction.

\begin{remark}
Suppose that we have done these operations and are in the situation of
bad reduction -- that is, with $g(0,0,0,1) \in \O^\times$ and $s \in
\{1,2\}$.  Then the singular point of $\Xs$ can never lift to a point
of $X(k)$.  For any such point would be of the form
$(x_0:x_1:x_2:x_3)$ with the $x_i \in \O$, $v(x_i) > 0$ for $i=0,1,2$
and $v(x_3)=0$.  But then $v(f(x_0,x_1,x_2)) \ge 3$, whereas $v(\pi^s
g(x_0,x_1,x_2,x_3))=s$, and so the equation~\eqref{eq:cone} cannot be
satisfied.
\end{remark}

Before we state the main theorem of this section, note the following
properties of plane sections of $X$.  If $H$ is a plane in $\P^3_k$,
then the reduction of $H$ is a plane $\mathcal{H}_s$ in $\P^3_\F$,
which may or may not pass through the point $(0:0:0:1)$.  If
$\mathcal{H}_s$ does pass through $(0:0:0:1)$, we will say that $H$ is
\emph{bad}; otherwise, $H$ is \emph{good}.  Let $H$ be a good plane;
the reduction of $C = X \cap H$ is contained in $(\Xs \cap
\mathcal{H}_s)$, which is a non-singular cubic curve.  Now $C$ has
dimension $1$; by~\cite[Corollary~9.10]{Hartshorne:AG}, the reduction
of $C$ also has dimension $1$, so is equal to $(\Xs \cap
\mathcal{H}_s)$; we deduce that $C$ is a non-singular plane cubic
curve with non-singular reduction.  Indeed, all such $C$ have
isomorphic reductions, since the non-singular plane sections of a cone
are all isomorphic.

\begin{theorem} \label{thm:cone}
Let $X$ be a smooth cubic surface defined over $k$.  Suppose that $X$
has a model $\X \subset \P^3_\O$ of the form~\eqref{eq:cone} with
$g(0,0,0,1)$ a unit in $\O$, and $3 \nmid s$.  Then
\begin{enumerate}
\item \label{planesec} $\H^0(k, \Pic\Xb) \cong \Z$, that is, the only
  divisor classes on $X$ defined over $k$ are multiples of the class
  of a plane section;
\item \label{h1} $\H^1(k, \Pic \Xb) \cong (\Br X / \Br k)$ is either
  trivial or isomorphic to $(\Z/3\Z)$ or $(\Z/3\Z)^2$;
\item \label{deg3} every element of $\Br X$ splits over the field
  $\knr(\pi^{1/3})$;
\item \label{inj} if $C$ is a good plane section of $X$, then the
  natural map $(\Br X / \Br k) \to (\Br C / \Br k)$ is injective.
\end{enumerate}
\end{theorem}

\begin{remark}
In part~(\ref{h1}) above, all three possibilities can occur.  Indeed,
both $(\Z/3\Z)$ and $(\Z/3\Z)^2$ can occur for diagonal cubic
surfaces~\cite[Proposition~1]{CTKS}.  To achieve $\Br X = \Br k$, a
careful inspection of the proof of Lemma~\ref{lem:reduction} below
shows that the following condition suffices: the Jacobian of the
reduction of the cubic curve $C$ has no non-trivial 3-torsion defined
over $\F$.
\end{remark}

Before proving Theorem~\ref{thm:cone}, we deduce further results.

\begin{proposition} \label{prop:surj}
Let $\al{A}_1, \dotsc, \al{A}_n \in \Br X$ be Azumaya algebras whose
images are linearly independent in the $\F_3$-vector space $\Br X /
\Br k$.  Fix a base point $Q \in X(k)$.  Then, for any $P \in X(k)$,
each $(\al{A}_i(P) - \al{A}_i(Q))$ lies in $(\Br k)[3]$, and the
product of the evaluation maps
\[
X(k) \to (\frac{1}{3}\Z / \Z)^n, \qquad P \mapsto (\inv \al{A}_1(P) -
\inv \al{A}_1(Q), \dotsc, \inv \al{A}_n(P) - \inv \al{A}_n(Q))
\]
is surjective.  In particular, if $\al{A}$ is any Azumaya algebra
on $X$ not equivalent to a constant algebra, then the map $\al{A}
\colon X(k) \to \Q / \Z$ takes three distinct values.
\end{proposition}
\begin{proof}
Firstly, note that the hypotheses and conclusion are unchanged if we
change any $\al{A}_i$ by a constant algebra; we will need to do this
below.

Let $C$ be a good plane section of $X$.  Note that $C(k) \neq
\emptyset$, since every smooth curve of genus 1 over a finite field
has a rational point, and this lifts by Hensel's Lemma.
Lichtenbaum~\cite{Lichtenbaum:IM-1969} has shown that the evaluation
pairing $C(k) \times \Br C \to \Br k$ extends to a non-degenerate
pairing $\Pic C \times \Br C \to \Br k$, and this gives rise to a
non-degenerate pairing
\begin{equation} \label{eq:licht2}
\Pic_0 C \times (\Br C / \Br k) \to \Br k \cong \Q/\Z \text{.}
\end{equation}
Since $C$ is an elliptic curve, we can identify $C(k)$ with $\Pic_0 C$
by choosing a base point $O \in C(k)$ and mapping $P \in C(k)$ to the
divisor class $[P-O] \in \Pic_0 C$.  We can thereby
view~\eqref{eq:licht2} as a pairing of Abelian groups
\[
C(k) \times (\Br C / \Br k) \to \Q/\Z
\]
defined as follows: given a class $\alpha$ in $\Br C / \Br k$, choose
the unique Azumaya algebra $\al{A}$ in that class satisfying
$\al{A}(O)=0$; then $(P, \alpha) = \al{A}(P)$.  For each integer $r$,
we obtain non-degenerate pairings of $(\Z/r\Z)$-modules
\[
C(k) / r C(k) \times (\Br C / \Br k)[r] \to (\frac{1}{r}\Z / \Z) \text{.}
\]

Now let $\al{A}_1, \dotsc, \al{A}_n$ be as above.  By part~(\ref{h1})
of Theorem~\ref{thm:cone}, each $\al{A}_i$ is of order 3 in $(\Br X /
\Br k)$; by part~(\ref{inj}) of Theorem~\ref{thm:cone}, the natural
map $\H^1(k,\Pic \Xb) \to \H^1(k, \Pic \Cb)[3]$ is an injective map of
$\F_3$-vector spaces, and so $\al{A}_1, \dotsc, \al{A}_n$ restrict to
Azumaya algebras $\al{A}'_1, \dotsc, \al{A}'_n$ on $C$ which are again
linearly independent in $(\Br C / \Br k)[3]$.  After possibly changing
each $\al{A}_i$ by a constant algebra (to ensure $\al{A}'_i(O)=0$),
we obtain evaluation maps which are $\F_3$-linear maps $C(k) / 3 C(k)
\to (\frac{1}{3}\Z / \Z)$ given by linearly independent elements of
the dual $\F_3$-vector space to $C(k) / 3 C(k)$.  The proof is
finished by the following easy lemma in linear algebra.
\end{proof}

\begin{lemma}
Let $K$ be any field, $V$ a finite-dimensional vector space over $K$,
and $\alpha_1, \dotsc, \alpha_n$ linearly independent elements of the
dual space $V^*$.  Then the linear map $\phi\colon V \to K^n$ defined
by $\phi(v) = ( \alpha_1(v), \dotsc, \alpha_n(v) )$ is surjective.
\end{lemma}
\begin{proof}
Extend $\alpha_1, \dotsc, \alpha_n$ to a basis $\alpha_1, \dotsc,
\alpha_m$ of $V^*$, and let $v_1, \dotsc, v_m$ be the dual basis in $V
= V^{**}$.  Given $(x_1, \dotsc, x_n) \in K^n$, we have
\[
\phi(x_1 v_1 + \dotsb + x_n v_n) = (x_1, \dotsc, x_n)
\]
and so $\phi$ is surjective.
\end{proof}

\begin{proposition}
Let $Y$ be a smooth cubic surface defined over a number field $L$,
such that the base change of $Y$ to one completion $L_v$ of $L$, with
$v \nmid 3$, satisfies the conditions of Theorem~\ref{thm:cone}.  Then
there is no Brauer--Manin obstruction to the existence of rational
points on $Y$.
\end{proposition}
\begin{proof}
Write $Y_v$ for the base change of $Y$ to $L_v$, let $\bar{L}_v$ be an
algebraic closure of $L_v$ and let $\bar{Y}_v$ denote the base change
of $Y$ to $\bar{L}_v$.  Let $G$ be the absolute Galois group of $L$;
the absolute Galois group $\Gal(\bar{L}_v/L_v)$ can be identified with
the decomposition group $G_v \subseteq G$.  Since $\Pic \bar{Y}_v$ is
generated by the 27 lines on $\bar{Y}_v$, all of which are defined
over an algebraic extension of $L$, there is a canonical isomorphism
$\Pic \bar{Y} \cong \Pic \bar{Y}_v$ which respects the Galois action.

Firstly, we show that the restriction map $\H^1(L, \Pic \bar{Y}) \to
\H^1(L_v, \Pic \bar{Y}_v)$ is injective.  Its kernel is $\H^1(G/G_v,
\Pic Y_{L_v})$.  But by part~(\ref{planesec}) of Theorem~\ref{thm:cone},
$\Pic Y_{L_v} \cong \Z$, with the trivial Galois action, and so
$\H^1(G/G_v, \Pic Y_{L_v})=0$.  We deduce that $\Br Y / \Br L$
injects into $\Br Y_v / \Br L_v$, and in particular is an
$\F_3$-vector space.

It follows that, for any $\al{A} \in \Br Y$ and any adelic point
$(P_w)_{w \in \places{L}} \in Y(\ad_L)$, the sum
\[
\sum_{w \in \places{L}} \inv_w \al{A}(P_w)
\] 
lies in $(\frac{1}{3}\Z / \Z)$.  For we know that $3\al{A}$ is
equivalent to a constant algebra $\al{B} \in \Br L$, and so
\[
3 \sum_{w \in \places{L}} \inv_w \al{A}(P_w)
= \sum_{w \in \places{L}} \inv_w (3 \al{A})(P_w)
= \sum_{w \in \places{L}} \inv_w \al{B} = 0 \text{.}
\]

If $\Br Y / \Br L = 0$, then every Azumaya algebra on $Y$ is
equivalent to a constant algebra, so there is no Brauer--Manin
obstruction and we are finished.  Otherwise, let $\al{A}_1, \dotsc,
\al{A}_n$ be a minimal set of generators for $\Br Y / \Br L$; then
they are linearly independent in $\Br Y / \Br L$ considered as a
vector space over $\F_3$.  

Let $(Q_w)_{w \in \places{L}} \in Y(\ad_L)$ be any adelic point of $Y$,
and define $x_1, \dotsc, x_n \in (\frac{1}{3}\Z / \Z)$ by
\[
x_i = \sum_{w \in \places{L}} \inv_w \al{A}_i (Q_w) \text{.}
\]

We now focus our attention on the place $v$.  By
Corollary~\ref{prop:surj}, the map
\begin{align*}
Y(L_v) &\to (\frac{1}{3}\Z / \Z)^n, \\
P_v &\mapsto (\inv_v \al{A}_1(P_v) - \inv_v \al{A}_1(Q_v), \dotsc, 
\inv_v \al{A}_n(P_v) - \inv_v \al{A}_n(Q_v))
\end{align*}
is surjective.  So we may choose a point $P_v \in Y(L_v)$ such that
\[
\inv_v \al{A}_i(P_v) - \inv_v \al{A}_i(Q_v) = -x_i
\]
for all $i$.  For all other places $w \in \places{L} \setminus \{ v
\}$, set $P_w = Q_w$.  We obtain
\begin{align*}
\sum_{w \in \places{L}} \inv_w \al{A}_i (P_w) 
&= \big( \sum_{w \in \places{L}} \inv_w \al{A}_i (Q_w) \big) + (\inv_v
\al{A}_i (P_v) - \inv_v \al{A}_i (Q_v)) \\
&= 0 \text{ for all $i$.}
\end{align*}
Now, since $\al{A}_1, \dotsc, \al{A}_n$ generate $\Br Y / \Br L$, we
have
\[
\sum_{w \in \places{L}} \inv_w \al{A} (P_w) = 0 \text{ for all }
\al{A} \in \Br Y
\]
and so there is no Brauer--Manin obstruction to the existence of a
rational point on $Y$.
\end{proof}

\begin{remark}
Suppose that $\Br Y / \Br L$ is not trivial, so that there exists an
Azumaya algebra $\al{A}$ on $Y$ which is not equivalent to a constant
algebra.  Then, in exactly the same way as above, we can construct an
adelic point $(P_w) \in Y(\ad_L)$ such that $\sum_{w \in \places{L}}
\inv_w \al{A} (P_w) \neq 0$, showing that there \emph{is} a
Brauer--Manin obstruction to weak approximation on $Y$.
\end{remark}

We will now begin the proof of Theorem~\ref{thm:cone}, by giving an
explicit description of the reductions of the 27 lines on $\Xb$ and
the Galois action on them.  This result is similar to Exercise~IV--80
of~\cite{EH:GS}, and I am grateful to Professor Harris for sketching
this proof.  Once again, $\Cs$ denotes the reduction of the plane
section $C \subset X$.

\begin{proposition} \label{prop:flexes}
Let $\X \subset \P^3_\O$ be a cubic surface of the
form~\eqref{eq:cone} such that $g(0,0,0,1)$ is a unit in $\O$, and $3
\nmid s$.  Let $C$ be a good plane section of $X$.  Then 
\begin{enumerate}
\item \label{defK} all 27 lines on $\Xb$ are defined over the tamely
  ramified cyclic extension $K=\knr(\pi^{1/3})$ of $\knr$;
\item \label{red} reduction takes the 27 lines on $\Xb$ three-to-one
  onto the nine lines in the cone $\Xbs$ lying over the nine points of
  inflection of the smooth plane cubic curve $\Cbs$;
\item \label{coplanar} each triple of lines on $\Xb$ with common
  reduction consists of three coplanar lines; and
\item \label{action} the Galois group $\Gal(K/\knr)$ acts cyclically
  on each triple of lines.
\end{enumerate}
\end{proposition}

\begin{proof}
The hypotheses and conclusions of the proposition are unchanged if we
replace $k$ by its maximal unramified extension $\knr$, and so we will
assume that $k=\knr$.  We may also assume that $C$ is the plane
section given by $X_3=0$, since any non-singular plane section of the
cone $\Xbs$ will have its flexes on the same lines.  Let $K=k(\Pi)$ be
the degree 3 ramified extension of $k$ with $\Pi^3=\pi$; the equation
of $\X$ becomes
\[
\X\colon f(X_0,X_1,X_2) + \Pi^{3s} g(X_0,X_1,X_2,X_3) = 0 \text{.}
\]
Although $X$ does not have good reduction, $X_K$ does; we will show
this by exhibiting an explicit change of variables producing a smooth
model of $X_K$.  Let $\phi\colon\P^3_K \to \P^3_K$ be the linear
automorphism given by the change of variables $X_i \mapsto \Pi^s X_i
\quad (i=0,1,2)$.  Applying $\phi$ to $X$ we obtain, as in the proof
of Lemma~\ref{lem:good}, a smooth model $\Y \subset \P^3_{\O_K}$ of
$X_K$ given by an equation of the form
\begin{equation} \label{eq:Y}
\mathcal{Y}\colon f(X_0,X_1,X_2) + aX_3^3 + \Pi h(X_0,X_1,X_2,X_3) = 0
\end{equation}
where $a = g(0,0,0,1) \in \O^\times$.  The generic fibre of $\Y$ is $Y
= \phi(X_K)$.

Since $\Y$ is smooth over $\O$, Lemma~\ref{lem:unram} shows that
$\Pic\bar{Y}$ is acted on trivially by $\Gal(\kb/k)$.  The 27 lines
on $\bar{Y}$ lie in distinct classes in $\Pic\bar{Y}$, so each is
defined over $k$, and so certainly over $K$.  As $\phi$ is defined
over $K$, we deduce that the 27 lines on $\Xb$ are defined over $K$,
                      proving part~(\ref{defK}).

To prove part~(\ref{red}), we will reduce it to the corresponding
statement for $Y$.  Notice that $\phi$ fixes the plane $H = \{ X_3 = 0
\} \subset \P^3_K$, so $\phi$ also fixes the curve $C_K = H \cap X_K$,
and each line $L \subset \Xb$ meets $C_K$ in the same point $P$ as its
image $\phi(L) \subset \bar{Y}$ does.  We would like to make the same
statement about reductions: the reduction of each line $L$ on $\Xb$
meets the reduction $\Cbs$ in the same point as the reduction of its
image $\phi(L) \subset \bar{Y}$ does.  The only problem is that, as
remarked in Section~\ref{sec:reduction}, the intersection of the
reductions of two varieties may be strictly larger than the reduction
of their intersection, and we need to rule out this possibility.  To
do this, note that, if a line and a smooth plane cubic curve are
contained in a cubic surface, singular or not, then they must meet
transversely in one point.  So the reduction of $L$ meets $\Cbs$
transversely in one point, which is the reduction of $P$; similarly,
the reduction of $\phi(L)$ meets $\Cbs$ transversely in the same
point.  To prove part~(\ref{red}), it remains to show that the
reductions of the 27 lines on $\bar{Y}$ all meet $\Cbs$ in its
inflection points, with three lines meeting at each inflection point.
This can be done by straightforward calculation, as follows.

Let $\omega$ denote a primitive third root of unity in $k$ and
$\tilde{\omega}$ its image in $\F$.  Looking at equation~\eqref{eq:Y},
we see that the special fibre $\Ys$ of $\Y$ is a triple cover of
$\P^2_\F=\{ X_3=0 \}$ branched over the smooth cubic curve $\{f=0\}$.
The tangent plane to $\Ybs$ at each inflection point of this curve
contains three straight lines, all passing through the inflection
point, and interchanged cyclically by the automorphism $X_3 \mapsto
\tilde{\omega} X_3$ of $\Ybs$.  These nine triples of lines are all
the 27 lines in $\Ybs$.  Since $\Y$ is smooth, the argument of
Lemma~\ref{lem:unram} shows that the reduction map $\Pic \bar{Y} \to
\Pic \Ybs$ is injective; since the 27 lines on $\bar{Y}$ lie in
distinct classes in $\Pic \bar{Y}$, their reductions are therefore
distinct lines in $\Ybs$.  So the 27 lines which we have found on
$\Ybs$ are precisely the reductions of the 27 lines in $\bar{Y}$.  We
deduce that the reductions of the 27 lines on $\Xb$ also meet the
curve $\{X_3=0\} \subset \Xs$ at its inflection points in threes,
proving part~(\ref{red}).

To see that the lines in each triple are coplanar, note that $\bar{Y}$
and $\Ys$, being non-singular cubic surfaces, both have exactly 45
sets of three coplanar lines.  The reductions of coplanar lines are
again coplanar; we deduce that three lines in $\bar{Y}$ are coplanar
if and only if their reductions in $\Ybs$ are coplanar.  We have seen
that the three lines in $\Ybs$ passing through any one inflection point
of $\Cbs$ are coplanar, and therefore so are the three lines of
$\bar{Y}$ reducing to them.  Since $\phi$ is a linear automorphism,
the corresponding triples of lines in $\Xb$ are also coplanar.  This
proves part~(\ref{coplanar}).

It remains to find the Galois action on the 27 lines of $\Xb$.  Let
$\sigma \in \Gal(K/k)$ be the automorphism such that $\sigma \Pi =
\omega \Pi$, and write $\sigma$ also for the corresponding
automorphism $\P^3_K \to \P^3_K$ of schemes over $k$.  Let $\psi\colon
\P^3_K \to \P^3_K$ be the linear automorphism $X_3 \mapsto \omega^s
X_3$.  Then a direct calculation shows that $\phi \sigma = \psi \sigma
\phi$.  Letting these automorphisms act on the sets of lines of $\Xb$
and $\bar{Y}$, we get a commutative diagram
\[
\begin{CD}
\{ \textrm{lines on } \Xb \} @>{\phi}>> \{ \textrm{lines on } \bar{Y}
\} \\
@V{\sigma}VV                            @V{\psi \sigma}VV \\
\{ \textrm{lines on } \Xb \} @>{\phi}>> \{ \textrm{lines on } \bar{Y}
\}
\end{CD} \text{.}
\]
As noted above, the lines on $\bar{Y}$ are all defined over $\knr$,
and so $\sigma$ acts trivially on them; we can rewrite $\psi\sigma$
as $\psi$ on the right-hand arrow.  Let $\tilde{\psi} \colon \P^3_\F
\to \P^3_\F$ be the linear automorphism $X_3 \mapsto \tilde{\omega}^s
X_3$; then we can extend our commutative diagram as follows:
\[
\begin{CD}
\{ \textrm{lines on } \Xb \} @>{\phi}>> \{ \textrm{lines on } \bar{Y}
\} @>{\textrm{reduce}}>> \{ \textrm{lines on } \Ybs \} \\
@V{\sigma}VV @V{\psi}VV @V{\tilde{\psi}}VV \\
\{ \textrm{lines on } \Xb \} @>{\phi}>> \{ \textrm{lines on } \bar{Y}
\} @>{\textrm{reduce}}>> \{ \textrm{lines on } \Ybs \} \\
\end{CD} \text{.}
\]
The horizontal maps are all bijections; so, to describe the action of
$\sigma$ on the lines of $\Xb$, we need to describe the action of
$\tilde{\psi}$ on the lines of $\Ybs$.  This was described above: given
that $3 \nmid s$, the action of $\tilde{\psi}$ is to cyclically
permute each triple of lines on $\Ybs$ passing through one inflection
point of the curve $\Cbs$.  We deduce that $\sigma$ cyclically
permutes each triple of lines on $\Xb$ with common reduction.
\end{proof}

\begin{remark}
The first part of the proposition can also be proven in a more
elementary way as follows.  Clebsch~\cite{Clebsch:JRAM-1861} gave a
covariant of degree nine associated to any cubic surface, which
vanishes precisely on the 27 lines in that surface.  Under the
hypotheses of the proposition, this covariant can be written as
\[
6 \pi^2 a^3 H_f^3 + (\text{multiple of }\pi^3)
\]
where $a = g(0,0,0,1)$ and $H_f$ is the Hessian of $f$, that is, the
determinant of the $3\times 3$ matrix of second-order partial derivatives
of $f$.  Since the Hessian of a smooth cubic curve vanishes precisely at
the nine flexes, this proves the result.  Unfortunately this approach
does not appear to give any easy way of showing that the triples of
lines are coplanar, nor of determining the Galois action.
\end{remark}

\begin{corollary} \label{cor:h0h1}
$\H^0(\knr, \Pic\Xb) \cong \Z$ and $\H^1(\knr, \Pic\Xb) \cong
  (\Z/3\Z)^2$.
\end{corollary}
\begin{proof}
To see that $\H^0(\knr, \Pic\Xb) \cong \Z$, observe that the only
Galois-fixed linear combinations of lines are made up of coplanar
triples.  But any coplanar triple of lines is a plane section of $X$,
and so they are all linearly equivalent.  The calculation of
$\H^1(\knr, \Pic\Xb)$ is Lemma~5 of~\cite{SD:MPCPS-1993}.
\end{proof}

We can now prove parts~(\ref{planesec})--(\ref{deg3}) of
Theorem~\ref{thm:cone}.

\begin{proof}[Proof of Theorem~\ref{thm:cone},
    parts~(\ref{planesec})--(\ref{deg3}).]  
By Corollary~\ref{cor:h0h1}, the only divisor classes defined over
$\knr$ are the multiples of the class of a plane section; so this is
certainly also true over $k$, proving part~(\ref{planesec}) of
Theorem~\ref{thm:cone}.

The inflation-restriction sequence
\[
\H^1(\knr/k, \H^0(\knr, \Pic\Xb)) \xrightarrow{\inf} \H^1(k, \Pic\Xb)
\xrightarrow{\res} \H^1(\knr, \Pic\Xb)
\]
shows that the restriction map is injective, since $\H^1(\knr/k,
\Z)=0$.  According to Corollary~\ref{cor:h0h1}, the group $\H^1(\knr,
\Pic\Xb)$ is isomorphic to $(\Z/3\Z)^2$, hence part~(\ref{h1}) of
Theorem~\ref{thm:cone}.

Finally, part~(\ref{deg3}) of Theorem~\ref{thm:cone} follows from the
fact that all the 27 lines on $\Xb$ are defined over $K$, the degree
$3$ extension of $\knr$.
\end{proof}

It remains to prove part~(\ref{inj}) of Theorem~\ref{thm:cone}, that
the natural map $(\Br X / \Br k) \to (\Br C / \Br k)$ is injective.

\begin{lemma} \label{lem:res}
To show that $(\Br X / \Br k) \to (\Br C / \Br k)$ is injective, it is
enough to show that the map $\H^1(\knr, \Pic\Xb) \to \H^1(\knr,
\Pic\Cb)$ is injective.
\end{lemma}
\begin{proof}
Since $X$ is rational, $\Br\Xb=0$; and, since $C$ is a curve,
$\Br\Cb=0$.  The exact sequence~\eqref{eq:hs} coming from the
Hochschild--Serre spectral sequence therefore gives a natural
isomorphism $\Br X/\Br k \cong \H^1(k,\Pic\Xb)$ and similarly for $C$.
So we must show that $\H^1(k, \Pic\Xb) \to \H^1(k, \Pic\Cb)$ is
injective.  The following diagram commutes:
\[
\begin{CD}
\H^1(k, \Pic\Xb)    @>>> \H^1(k, \Pic\Cb) \\
@V{\res}VV              @VV{\res}V \\
\H^1(\knr, \Pic\Xb) @>>> \H^1(\knr, \Pic\Cb)
\end{CD} \text{.}
\]
We showed in the proof of part~(\ref{h1}) of Theorem~\ref{thm:cone}
that the left-hand map is injective.  If the bottom map is also
injective, then their composition is injective and it follows that the
top map is injective.
\end{proof}

Let $\Cs$ be the reduction of the curve $C$; since $C$ is a good plane
section of $X$, then $\Cs$ is a plane section of the cone $\Xs$ which
is a smooth plane cubic curve.  There is a reduction map $\Pic \Cb \to
\Pic \Cbs$.  To prove our desired result that the map $\H^1(\knr, \Pic
\Xb) \to \H^1(\knr, \Pic \Cb)$ is injective, we will prove a statement
which at first sight looks stronger.

\begin{lemma} \label{lem:reduction}
The composition
\[
\H^1(\knr, \Pic\Xb) \to \H^1(\knr, \Pic\Cb) \to \H^1(\knr, \Pic \Cbs)
\]
is injective (where $\Gal(\kb/\knr)$ acts trivially on $\H^1(\knr,
\Pic \Cbs)$).
\end{lemma}

\begin{proof}
We first re-interpret the composition of maps $\phi\colon \Pic \Xb \to
\Pic \Cb \to \Pic \Cbs$.  By~\cite[Proposition~20.3]{Fulton:IT}, the
following diagram commutes:
\[
\begin{CD}
\Pic \Xb @>{\text{reduce}}>> \Cl \Xbs \\
@VVV                         @VVV \\
\Pic \Cb @>>{\text{reduce}}> \Pic \Cbs
\end{CD}
\]
and so $\phi$ can also be obtained by first taking the reduction of a
divisor on $\Xb$ and then restricting to $\Cbs$.  Since $\Pic \Xb$ is
generated by the 27 lines, this is precisely what
Proposition~\ref{prop:flexes} tells us about.

Let $\phi_*$ denote the induced map on cohomology groups; we shall
show that $\phi_*$ identifies $\H^1(\knr, \Pic\Xb)$ with the 3-torsion
subgroup of $\H^1(\knr, \Pic\Cbs)$.  Let $K/\knr$ be the extension
described in Proposition~\ref{prop:flexes}; then $K$ is the unique
cyclic extension of $\knr$ of degree $3$.  There is a commutative
diagram
\begin{equation} \label{eq:inf}
\begin{CD}
\H^1(\knr, \Pic\Xb)    @>{\phi_*}>> \H^1(\knr, \Pic \Cbs) [3] \\
@A{\inf}AA                        @AA{\inf}A \\
\H^1(K/\knr, \Pic X_K) @>{\phi_*}>> \H^1(K/\knr, \Pic \Cbs)
\end{CD}
\end{equation}
in which both inflation maps are isomorphisms, for the following
                               reasons:
\begin{itemize}
\item The 27 lines on $\Xb$ are defined over $K$, so $\Gal(\kb/K)$
  acts trivially on $\Pic\Xb$ and therefore $\H^1(K, \Pic \Xb)=0$; the
  inflation-restriction sequence then shows that the left-hand
  inflation map is an isomorphism.
\item The right-hand inflation map can be identified with the natural
  homomorphism
\[
\Hom(\Gal(K/\knr), \Pic \Cbs) \to \Hom(\Gal(\kb/\knr), \Pic \Cbs)[3]
\]
which is an isomorphism because $\Gal(K/\knr)$ is the maximal quotient
of $\Gal(\kb/\knr)$ killed by 3.
\end{itemize}
It remains to prove that the bottom horizontal arrow in~\eqref{eq:inf}
is an isomorphism.  We already know that $\H^1(K/\knr, \Pic X_K) \cong
(\Z/3\Z)^2$, from Corollary~\ref{cor:h0h1}.  Also, $\H^1(K/\knr, \Pic
\Cbs)$ is naturally isomorphic to the 3-torsion in $\Pic \Cbs$, which
is isomorphic to $(\Z/3\Z)^2$.  So, to prove that this map is an
isomorphism, it is enough to show that it is surjective.

Let $M$ denote the kernel of $\phi$ and $N$ the image.  Then
Proposition~\ref{prop:flexes} shows that $N$ is the subgroup of $\Pic
\Cbs$ generated by the inflection points.  It is well known that the
differences of the inflection points generate the 3-torsion of
$\Pic\Cbs$; in particular, $N$ contains all the 3-torsion in $\Pic
\Cbs$, and so $\H^1(K/\knr, N) = \H^1(K/\knr, \Pic \Cbs)$.  We have
\[
\begin{CD}
\H^1(K/\knr, \Pic\Xb) @>>> \H^1(K/\knr, N) @>>>
\H^2(K/\knr, M) \\
@. @| \\
@. \H^1(K/\knr, \Pic \Cbs)
\end{CD} \text{.}
\]
We will show that $\H^2(K/\knr, M)=0$.  By Corollary~\ref{cor:h0h1},
$\H^0(K/\knr, \Pic\Xb)$ is generated by the class $H$ of a plane
section.  For any non-zero integer $n$, $\phi(nH)$ is a divisor on
$\Cbs$ of non-zero degree, so certainly not principal; we deduce that
$\H^0(K/\knr, \Pic\Xb)$ injects into $\H^0(K/\knr, N)$.  Looking at
the exact sequence
\[
0 \to \H^0(K/\knr, M) \to \H^0(K/\knr, \Pic\Xb) \to \H^0(K/\knr, N)
\]
we see that $\H^0(K/\knr,M)=0$.  Since $\Gal(K/\knr)$ is cyclic, we
have
\[
\H^2(K/\knr, M) \cong \hat{\H}^0(K/\knr, M) = \H^0(K/\knr,M) /
\mathrm{N}_{K/\knr} M
\]
and so $\H^2(K/\knr,M)=0$, completing the argument.
\end{proof}

\begin{remark}
If the special fibre of $\X$ is a cone, but the total space $\X$ is
too singular for Proposition~\ref{prop:flexes} to hold, then the
reduction map $\phi\colon \Pic\Xb \to \Pic\Cbs$ can indeed behave
differently.  We give an example.  Let $\X$ be the cubic surface over
$\Z_5$ defined by the equation
\[
X_0^3 + X_1^3 + X_2^3 + 5X_2 X_3^2 = 0 \text{.}
\]
This is of the form~\eqref{eq:cone}, but with $g(0,0,0,1)=0$.  Using a
computer algebra system, one can compute the Fano scheme of lines on
$X$ over $\Q_5$, take its flat completion over $\Z_5$ and look at the
resulting special fibre; this describes the reductions of the 27 lines
on $\Xb$.  We find that the reductions are supported on 12 lines in
the cone $\Xbs$, and the associated points of the curve $\Cbs$
constitute a subgroup isomorphic to $\Z/2\Z \times \Z/6\Z$ in its
Jacobian.
\end{remark}

\subsection*{Acknowledgements}
This research was carried out during the 2007--8 Warwick EPSRC
Symposium on Algebraic Geometry.  I thank the Mathematics Research
Centre of the University of Warwick, and in particular Samir Siksek,
Ronald van Luijk and Miles Reid, for their hospitality.  I also thank
both Jean-Louis Colliot-Th\'el\`ene and the anonymous referees for
suggesting several improvements to this article.

\bibliographystyle{abbrv} 
\bibliography{martin}

\end{document}